%% file: main.tex
\documentclass[11pt]{article}
\input{header}

\input{marco}

\title{Gradient Descent Achieves $\mathcal{O}(\eps^{-7/4})$ Convergence}
\author{Chang He\thanks{Department of Industrial and Systems Engineering, University of Minnesota. \texttt{hec@umn.edu}}\and Xiaotian Jiang\thanks{Department of Industrial and Systems Engineering, University of Minnesota. \texttt{jian0851@umn.edu}}}
\date{}

\begin{document}
\maketitle

\begin{abstract}
We propose a deterministic stepsize rule for gradient descent on a lower-bounded nonconvex function with Lipschitz-continuous gradient and Hessian. Classical analysis shows that a fixed stepsize guarantees a decrease of order $\eps^2/L$ per update. In contrast, our rule uses a sequence of $\cO(m)$ updates with a total stepsize of order $m^2/L$. We establish a weighted descent inequality for every prefix of the update sequence on a quadratic objective and control the error caused by Hessian variation on a general objective. This error determines the choice $m\asymp(L^2/(\rho\eps))^{1/4}$ for small $\eps$. We show that gradient descent with our stepsize rule finds an $\eps$-stationary point within $\cO(\eps^{-7/4})$ gradient evaluations, matching the optimal gradient-complexity guarantee.

\end{abstract}

\section{Introduction}
We consider unconstrained minimization of a lower-bounded, possibly nonconvex function $f:\bR^d\to\bR$. Given $\eps>0$, the goal is to find an iterate satisfying $\norm{\nabla f(x)}\le\eps$. When $\nabla f$ is $L$-Lipschitz, gradient descent with stepsize $1/L$ takes $\cO(\eps^{-2})$ gradient evaluations \citep{nesterov2018lectures}. If the Hessian is also $\rho$-Lipschitz, \cite{carmon2017} obtain $\cO(\eps^{-7/4}\log(1/\eps))$ with accelerated iterations and negative-curvature exploitation. Restarted accelerated gradient methods attain $\cO(\eps^{-7/4})$ \citep{lilin2023,marumo2024}. We ask whether choosing the scalar stepsize alone can attain this exponent.

With stepsize $1/L$, the descent guaranteed by one update is proportional to $\norm{\nabla f(x)}^2/L$. Increasing a single stepsize can cause ascent even on a one-dimensional quadratic. We instead bound the total decrease over an epoch. A sequence of $\cO(m)$ updates whose coefficients sum to order $m^2$ can produce a decrease of order $\eps^2m^2/L$ before a small gradient is found. If an epoch realizes this decrease, its average decrease per update improves by a factor of $m$. The order of the coefficients matters because a large step can amplify errors from earlier updates. On quadratic objectives, this dependence on the ordering also appears in the analysis of fractal stepsize schedules \citep{agarwal2021fractal}.

For a general objective, the main issue is the change in the Hessian during the epoch. We insert shorter sequences around the largest coefficients and prove that the resulting table amplifies gradient errors by at most $\cO(m^2)$. Clipping the stepsizes using observed gradients bounds the path length of each epoch with a long step by $\cO(\eps m^2/L)$. Thus $\nabla^2 f/L$ can change by at most $\cO(\rho\eps m^2/L^2)$. Together these bounds yield the dimensionless error term $\rho\eps m^4/L^2$. Choosing $m$ of order $(L^2/(\rho\eps))^{1/4}$ keeps this term bounded and gives the $\eps^{-7/4}$ complexity. The proof accounts for every gradient evaluation, including those in epochs whose descent test fails.

\section{Preliminaries}\label{sec:preliminaries}
Let $f:\bR^d\to\bR$ be twice continuously differentiable and bounded below, and set $\Delta=f(x^0)-\inf f$. The notation $\norm{\cdot}$ denotes the Euclidean norm of a vector and the induced operator norm of a matrix. All constants implicit in $\cO(\cdot)$, $\lesssim$ and $\asymp$ are absolute. We take an empty sum to be zero and an empty product to be one.
\begin{definition}\label{def:stationary}
An point $x$ is $\eps$-stationary if $\norm{\nabla f(x)}\le\eps$.
\end{definition}
\begin{assumption}\label{ass:gradient}

There exists $L>0$ satisfying $\norm{\nabla f(x)-\nabla f(y)}\le L\norm{x-y}$ for all $x,y\in\bR^d$.

\end{assumption}
\begin{assumption}\label{ass:hessian}
There exists $\rho>0$ satisfying $\norm{\nabla^2 f(x)-\nabla^2 f(y)}\le\rho\norm{x-y}$ for all $x,y\in\bR^d$.
\end{assumption}

To see where the exponent $-2$ comes from, write a gradient descent step as $x^+=x-(h/L)\nabla f(x)$ with $h>0$. Under \Cref{ass:gradient}, the smoothness inequality gives
\begin{equation}\label{eq:fixed}\textstyle f(x)-f(x^+)\ge(h/L)(1-h/2)\norm{\nabla f(x)}^2.\end{equation}
For $h=1$, each update before termination decreases $f$ by more than $\eps^2/(2L)$. Summation gives at most $2L\Delta/\eps^2$ updates. For $f(x)=Lx^2/2$, equality holds in \eqref{eq:fixed}, and every $h>2$ increases $f$ at $x\ne0$.

To improve the count, suppose an epoch takes at most $\cO(m)$ updates and produces a decrease of order $\eps^2m^2/L$. Dividing the available decrease by the decrease per epoch and counting one final epoch gives the target
\begin{equation}\label{eq:target}\textstyle \underbrace{\cO(m)}_{\text{updates per epoch}}\,\underbrace{\cO\bigl(1+L\Delta/(\eps^2m^2)\bigr)}_{\text{number of epochs}}=\cO\bigl(m+L\Delta/(\eps^2m)\bigr).\end{equation}
Section~\ref{sec:quadratic} derives a scalar condition under which the total stepsize produces this decrease on a quadratic. Section~\ref{sec:nonlinear} gives the clipping rule and test that make the same counting argument possible for a varying Hessian.

\section{Quadratic Acceleration}\label{sec:quadratic}
We seek a schedule whose total coefficient sum is of order $m^2$ in $\cO(m)$ updates and whose first $t$ steps guarantee sufficient decrease for every $t$. To determine what this requires, take any quadratic $Q$ with $\norm{\nabla^2Q}\le L$ and write $x_{i+1}=x_i-(h_i/L)g_i$, where $g_i=\nabla Q(x_i)$ and $h_i>0$. For $p_i(z)=\prod_{j<i}(1-h_jz)$, we claim that the scalar inequality on the left guarantees the corresponding decrease on the right,
\begin{equation}\label{eq:quadratic-prefix}\textstyle z\sum_{i<t}h_i p_i(z)^2+p_t(z)^2\le1\ \text{for all }z\in[0,1]\ \Longrightarrow\ Q(x_0)-Q(x_t)\ge\tfrac1{2L}\sum_{i<t}h_i\norm{g_i}^2.\end{equation}

To prove the decrease in \eqref{eq:quadratic-prefix}, consider one eigendirection of $\nabla^2Q/L$ with eigenvalue $z$. Its gradient coordinate is $g_i(z)=p_i(z)g_0(z)$, and its contribution to the decrease is $L^{-1}\sum_{i<t}h_i(1-zh_i/2)g_i(z)^2$ by exact quadratic expansion. For $z\le0$, this is at least $L^{-1}\sum_{i<t}h_i g_i(z)^2$. For $z>0$, the identity $p_i^2-p_{i+1}^2=2zh_i(1-zh_i/2)p_i^2$ makes the contribution equal to $g_0(z)^2(1-p_t(z)^2)/(2Lz)$. The scalar condition bounds this below by $(2L)^{-1}\sum_{i<t}h_i g_i(z)^2$. Summing over an orthonormal eigenbasis proves \eqref{eq:quadratic-prefix}. If all gradients exceed $\eps$, a full table therefore gives decrease greater than $(\eps^2/(2L))\sum_i h_i$. We next construct coefficients with the required sum and prefix bound.

\paragraph{Recursive construction.} On a quadratic, a step with coefficient $a$ multiplies the gradient in a direction of normalized curvature $z$ by $1-az$. To split a step into two, require $(1-bz)(1-cz)=1-4az(1-z)$. Matching coefficients gives $b+c=bc=4a$. We put the smaller root first. Starting from $a_0=1$ and $a_1=2$ gives the raw table $A_n=(a_0,\ldots,a_{n-1})$ with
\begin{equation}\label{eq:raw}\textstyle a_{2j}=2/(1+\sqrt{1-1/a_j}),\ a_{2j+1}=4a_j-a_{2j},\ j\ge1.\end{equation}
Splitting a whole layer doubles its length and quadruples its sum. The first $n$ raw coefficients therefore have sum of order $n^2$, as required by \eqref{eq:target}. The issue for a nonlinear function is that an error introduced partway through a long step sequence sees only the remaining factors. We limit that amplification by inserting a smaller copy of the table on each side of the large raw coefficients.

Let $n$ be a power of two. For $j=0,\ldots,n-1$, let $q_j$ be the largest power of two dividing $j+1$. Set $r_j=q_j^2/(4n)$ if $q_j^2\ge4n$, and $r_j=0$ otherwise. To construct $F_n$, read the raw coefficients $a_0,\ldots,a_{n-1}$ in order. Immediately before and after each $a_j$, insert the entries of $F_{r_j}$ in their original order. This gives the ordered list
\begin{equation}\label{eq:table}\textstyle F_n=\bigl(F_{r_0},a_0,F_{r_0},\ldots,F_{r_{n-1}},a_{n-1},F_{r_{n-1}}\bigr).\end{equation}
Here $F_0$ is empty, so $r_j=0$ contributes only $a_j$. Every positive $r_j$ is smaller than $n$, so its table has already been constructed. For example, $F_1=(1)$ and $F_2=(1,2)$. At $n=4$, $(r_0,r_1,r_2,r_3)=(0,0,0,1)$, so only $a_3$ receives insertions and $F_4=(1,2,4-2\sqrt2,1,4+2\sqrt2,1)$. To see why $r_j$ scales with $q_j^2$, consider a raw interval whose largest scale is $q_j$ in a direction of normalized curvature $0\le z\le1$. Its amplification is $\cO(1+zq_j^2)$, and a complete adjacent table contributes a factor of order $1/\sqrt{1+zr_j^2}$ or smaller. For $r_j>0$, the identity $q_j^2=4nr_j$ bounds their product by $\cO(1+n\sqrt z)$. Appendix~\ref{sec:kernel} sums these effects over all possible error locations.

\begin{figure}[htbp]
\centering
\begin{tikzpicture}
\begin{semilogyaxis}[
  scale only axis, width=.84\linewidth, height=5.8cm,
  xmin=-.75, xmax=29.85, ymin=.78, ymax=141, log origin y=0,
  xtick={0,4,8,12,16,20,24,29}, ytick={1,2,5,10,25,50,100},
  yticklabels={1,2,5,10,25,50,100}, minor y tick num=0,
  xlabel={Step index $i$}, ylabel={Normalized stepsize $b_i$},
  axis lines=left, axis line style={black!70,-},
  tick align=outside, tick style={black!55}, major tick length=2.5pt,
  tick label style={font=\small}, label style={font=\small},
  ymajorgrids, grid style={draw=black!13,line width=.35pt}, enlargelimits=false]
\addplot+[ycomb,EvenRed,line width=.85pt,mark=square*,mark size=2.25pt,
  mark options={solid,fill=EvenRed,draw=EvenRed}]
  coordinates {(0,1)(1,2)(2,1.17157288)(3,6.82842712)(4,1.44646269)(5,3.23982881)
    (6,1.03956613)(7,1)(8,26.27414237)(9,1)(10,1.28570215)(11,4.50014861)
    (12,1.09201921)(13,11.86729602)(14,1.67351368)(15,2.48475084)(16,1.00970056)
    (17,1)(18,2)(19,1.17157288)(20,1)(21,6.82842712)(22,1)(23,104.08686892)
    (24,1)(25,2)(26,1.17157288)(27,1)(28,6.82842712)(29,1)};
\node[anchor=west,text=EvenRed,font=\small,inner sep=3pt] at (axis cs:23.25,104.08686892) {$104.09$};
\end{semilogyaxis}
\end{tikzpicture}
\caption{The quadratic stepsize schedule $F_{16}=(b_0,\ldots,b_{29})$ from \eqref{eq:table}. The squares show the normalized stepsizes $b_i$ on a logarithmic scale. The actual stepsizes are $b_i/L$. The copies of $F_4$ at indices $17$--$22$ and $24$--$29$ flank the longest step at $i=23$.}\label{fig:quadratic-steps}
\end{figure}
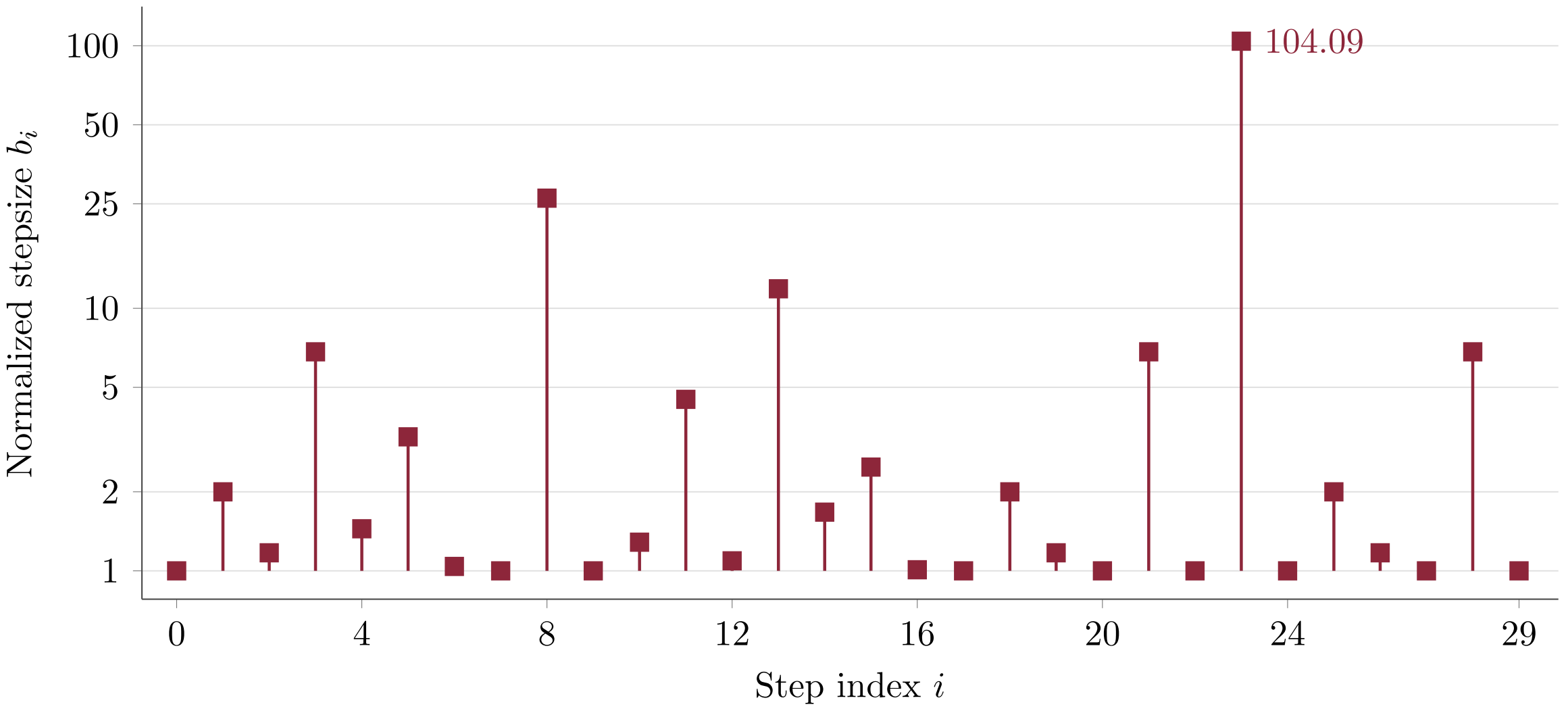

Figure~\ref{fig:quadratic-steps} shows how the insertions change the schedule at $n=16$. The 16 raw coefficients sum to 171. Inserting two copies each of $F_1$ and $F_4$ gives 30 steps with total stepsize $199/L$.

\FloatBarrier

\section{From Quadratic to Nonlinear Functions}\label{sec:nonlinear}
To control the change in the Hessian, each epoch limits how far the iterates can travel. Write the stepsize table as $F_m=(b_0,\ldots,b_{J-1})$, index the iterates from $x_0$ within each epoch, and write $g_i=\nabla f(x_i)$. Set the threshold $B=(\eps^2/2)\sum_{j=0}^{J-1}b_j$ for the first epoch. Each step contributes $h_i\norm{g_i}^2$ to the sum compared with $B$. The epoch always starts with the ordinary step $x_1=x_0-g_0/L$. If $\norm{g_0}^2\ge B$, this step ends the epoch. Otherwise, each subsequent coefficient is clipped so that the sum stays at most $B$. Thus
\begin{equation}\label{eq:update}\textstyle h_0=1,\ h_i=\min\bigl\{b_i,\tfrac{B-\sum_{j<i}h_j\norm{g_j}^2}{\norm{g_i}^2}\bigr\}\ (i\ge1),\ x_{i+1}=x_i-(h_i/L)g_i.\end{equation}
At the end of an epoch, the gradients at the two ends of each step estimate the decrease along that step by the trapezoidal rule. The algorithm halves $B$ for the next epoch if the sum of these estimates is less than $(1/(4L))\sum_i h_i\norm{g_i}^2$ over the completed epoch. The test below writes this comparison directly in terms of the observed gradients.

\begin{algorithm}[htbp]
\caption{Gradient descent with adaptive stepsizes}\label{alg:gd}
\begin{algorithmic}[1]
\Require Starting point $x^0$, smoothness bounds $L,\rho>0$, tolerance $\eps>0$.
\State Let $m$ be the largest power of two at most $\max\{1,(L^2/(\rho\eps))^{1/4}\}$ and form $F_m=(b_0,\ldots,b_{J-1})$ by \eqref{eq:table}.
\State Set $B=(\eps^2/2)\sum_{j=0}^{J-1}b_j$. Set $x_0=x^0$ and evaluate $g_0=\nabla f(x_0)$.
\State \textbf{if} $\norm{g_0}\le\eps$ \textbf{then return} $x_0$.
\Loop
    \For{$i=0,\ldots,J-1$}
        \State Use $h_0=1$ for the first step. For $i\ge1$, use $h_i=\min\{b_i,(B-\sum_{j<i}h_j\norm{g_j}^2)/\norm{g_i}^2\}$.
        \State Set $x_{i+1}=x_i-(h_i/L)g_i$ and evaluate $g_{i+1}=\nabla f(x_{i+1})$.
        \State \textbf{if} $\norm{g_{i+1}}\le\eps$ \textbf{then return} $x_{i+1}$.
        \If{$\sum_{j=0}^{i}h_j\norm{g_j}^2\ge B$}
            \State Halve $B$ if $\sum_{j=0}^{i}h_j\inner{g_j}{g_j+2g_{j+1}}<0$.
            \State Set $(x_0,g_0)=(x_{i+1},g_{i+1})$ and \textbf{break} to start the next epoch.
        \EndIf
    \EndFor
\EndLoop
\end{algorithmic}
\end{algorithm}
\FloatBarrier

\paragraph{Descent and gradient complexity.} To obtain \eqref{eq:target}, we need a bound on the decrease in passing epochs, the possible increase in failed epochs, and how often $B$ is halved. The following proposition supplies all three.

\begin{proposition}\label{prop:epoch}
Suppose $m\ge2$. Let $B_0=(\eps^2/2)\sum_{j=0}^{J-1}b_j$ denote the initial value of $B$. Every completed epoch that passes the test decreases $f$ by at least $3B/(16L)$. On an epoch that fails, the increase in $f$ is at most a constant times $B/L$. There is an absolute $c\in(0,1)$ such that every completed epoch with $B\le cB_0$ passes.
\end{proposition}

The proof in Appendix~\ref{app:descent} bounds the error of the trapezoidal estimate and shows that, once $B$ is sufficiently small, the frozen-Hessian descent dominates the accumulated gradient errors. Each failure halves $B$. Hence \Cref{prop:epoch} bounds the number of failures by an absolute constant and keeps $B$ above a fixed fraction of $B_0$. The remaining epochs each decrease $f$ by order $\eps^2m^2/L$ at a cost of at most $4m$ updates. Counting them and the terminal epoch gives the following guarantee.

\begin{theorem}\label{thm:main}
Under \Cref{ass:gradient,ass:hessian}, for every $\eps>0$, Algorithm~\ref{alg:gd} returns an iterate $x$ with $\norm{\nabla f(x)}\le\eps$ after at most
\begin{equation}\label{eq:rate}\textstyle \cO\bigl(1+(L^2/(\rho\eps))^{1/4}+\Delta L^{1/2}\rho^{1/4}\eps^{-7/4}\bigr)\end{equation}
exact gradient evaluations.
\end{theorem}

Appendix~\ref{app:complexity} proves the evaluation bound for every tolerance, including the case $m=1$ where the method uses the ordinary stepsize $1/L$.

\FloatBarrier
\section{Conclusion}\label{sec:conclusion}
We developed this result through several rounds of interaction with GPT-6.2 Astra Ultra while refining the stepsize rule and its proof. Although the construction is algebraically involved, its main idea is to exploit Hessian Lipschitz continuity by viewing the objective function locally as a quadratic function with a cubic error term. The design proceeds in two steps. We first choose a sequence of scalar stepsizes that accelerates gradient descent on quadratic objectives. For a general nonlinear objective, we then clip each step using previously computed gradients. Clipping limits the distance traveled within an epoch. Hessian Lipschitz continuity controls the error in the quadratic model, ensuring that each completed epoch decreases the objective when $B$ is sufficiently small. We believe that the stepsize rule can be simplified further. The accompanying Lean~4 verification is available at \url{https://chang6234.github.io/blog.html}.

\bibliographystyle{abbrvnat}
\bibliography{ref}
\begin{appendices}
\numberwithin{equation}{section}
\section{Analysis of the Stepsize Table}\label{app:scalar}

To prove the epoch count we need bounds on the length and sum of the table. To prove descent we also need a prefix estimate for the frozen gradients and bounds on the products connecting two errors. For a list $W=(h_0,\ldots,h_{s-1})$ write $N(W)=s$ and $M(W)=\sum_{i<s}h_i$. Fix $0\le z\le1$, write $p_W=\prod_{i<s}(1-zh_i)$ and $E_W=z\sum_{i<s}h_i p_{W_{<i}}^2$, and take an empty product to be one. Here $W_{<i}$ and $W_{>i}$ are the portions strictly before and after index $i$.

\subsection{Raw coefficients and interval products}\label{app:raw}
\begin{lemma}\label{lem:raw}
For the raw table $A_n=(a_0,\ldots,a_{n-1})$ of \eqref{eq:raw}, $M(A_n)=(2n^2+1)/3$, $a_j\le 2q_j^2$, and $\max_{j<n}a_j\le n^2/2$ for $n\ge2$. For $0\le z\le1$, every prefix $W$ of $A_n$ satisfies $E_W(z)+p_W(z)^2\le1$. The complete product satisfies $|p_{A_n}(z)|\le(2n\sqrt z)^{-1}$ for $z>0$. For every nonempty consecutive raw interval $I$, $|p_I(z)|\le 1+2z\max_{j\in I}q_j^2$.
\end{lemma}
\begin{proof}
We need the sum and size of the coefficients to count updates, the prefix inequality for frozen-Hessian descent, and the two product bounds to control perturbations. Each estimate follows from the same splitting identity. If $b,c$ are the two children of $a$ in \eqref{eq:raw}, then $1\le b\le2\le c$, $b+c=bc=4a$, and $(1-bz)(1-cz)=1-ay$ with $y=4z(1-z)$. For the mass, each generation multiplies the sum by four, so the dyadic layer $[d,2d)$ has mass $2d^2$. Summing the layers and the initial unit gives $(2n^2+1)/3$. For the coefficient bounds, a left child is at most two and a right child is at most four times its parent. The number of consecutive right branches since the last left branch or the root $a_1=2$ is at most the number $\log_2 q_j$ of trailing ones in the binary representation of $j$. For $1\le j<n$ the total depth below $a_1$ is at most $\log_2 n-1$. These two counts give the stated bounds.

To prove the prefix inequality, remove the initial unit and first prove $(1-z)E_W(z)+p_W(z)^2\le1$ for every prefix of the remaining list. Passing from $n$ to $2n$ prepends $2$ to the child pairs of that list. For a prefix ending after the pairs of an old prefix $W$, we can apply induction at $y$ because $0\le y\le1$. In fact $b+c(1-bz)^2=4a(1-2z+bz^2)\le4a(1-y/2)$, so the new expression is at most
\[\textstyle \tfrac y2+(1-y)\bigl[(1-y/2)E_W(y)+p_W(y)^2\bigr]\le\tfrac y2+(1-y/2)\bigl[(1-y)E_W(y)+p_W(y)^2\bigr]\le1.\]
If the prefix also includes the next small child, the induction bound persists because $bz(1-z)+(1-bz)^2\le1$. The empty prefix starts the induction. Restoring the leading unit gives $E_{(1,W)}+p_{(1,W)}^2=z+(1-z)^2(E_W+p_W^2)\le z+(1-z)\bigl[(1-z)E_W+p_W^2\bigr]\le1$.

For the complete product, put $z=\sin^2 u$ with $0\le u\le\pi/2$. Pair compression replaces $\sin^2 u$ by $\sin^2(2u)$, so the product on the layer $[d,2d)$ is $\cos(2du)$. The sine doubling identity now gives $p_{A_n}(z)=\cos u\,\sin(2nu)/(2n\sin u)$ for $z>0$, which proves the claim.

To bound $p_I$, choose the unique $j\in I$ with largest $q_j=q$. Two indices of the same scale have a larger scale between them. The indices in $I$ through $j$ lie in $[j-q+1,j]$, and those after $j$ form a proper prefix of $[j+1,j+q]$. It therefore suffices to bound products on a block of $d$ descendants of a single coefficient $a$. We claim that every proper prefix has modulus at most one and every interval has modulus at most $\max(1,ad^2z)$. Repeated pair compression gives the full product $1-aw$ for some $0\le w\le\min(1,d^2z)$, since each compression maps $w$ to $4w(1-w)$. A proper prefix entering the right half contains a complete left half, whose root $b\le2$ gives modulus at most one, followed by a proper prefix. This proves the first claim by induction on $d$.

For the interval claim, the singleton case follows from $|1-az|\le\max(1,az)$. An interval in a single half obeys the induction bound since both child roots are at most $4a$. For a crossing interval, the part in the left half has modulus at most $\max(1,bd^2z/4)$. If the part in the right half is a proper prefix, its modulus is at most one. If it fills the right half, its modulus is at most both $c-1$ and $\max(1,cd^2z/4)$ by the full product formula. For $bd^2z/4\le1$ the latter bound suffices. For $bd^2z/4>1$ the product is at most $b(c-1)d^2z/4=cd^2z/4$. This proves the second claim.

Returning to $I$, the block ending at $j$ starts at an even multiple of $q$. When the left endpoint of this block is positive, its root has even index, hence coefficient at most two. When it starts at zero, it consists of the initial unit and complete dyadic layers rooted at $2$. Only the layer in which $I$ begins can amplify the product, and that layer has at most $q$ leaves. The part of $I$ after $j$ is a proper prefix of a single block and has modulus at most one. Multiplication gives $|p_I|\le\max(1,2q^2z)\le1+2q^2z$, as desired.
\end{proof}

\subsection{Stability under recursive insertion}\label{sec:stability}
\paragraph{Length and mass.} We claim $N(F_n)\le4n$ and $(2n^2+1)/3\le M(F_n)\le2n^2$. By \eqref{eq:table}, it suffices to control the total length and mass of the inserted copies. The scale $q=n$ occurs once and $q=n/2^k$ occurs $2^{k-1}$ times for $1\le k\le\log_2 n$. Since $r_j\le q_j^2/(4n)$, summing by scale gives
\[\textstyle 2\sum_{j<n}r_j\le\tfrac n2\bigl(1+\tfrac12\sum_{k\ge1}2^{-k}\bigr)=\tfrac34 n,\ 2\sum_{j<n}r_j^2\le\tfrac{n^2}8\bigl(1+\tfrac12\sum_{k\ge1}2^{-3k}\bigr)=\tfrac{15}{112}n^2 .\]
Thus induction gives $N(F_n)\le n+4(3n/4)=4n$ and $M(F_n)\le M(A_n)+2(15n^2/112)\le2n^2$. The lower bound is $M(A_n)$ by \Cref{lem:raw}. The same lemma bounds every coefficient of $F_n$ by $n^2/2$ for $n\ge2$, since inserted tables have smaller degree and $F_1=(1)$.

\paragraph{Stability under insertion.} We need the raw prefix bound of \Cref{lem:raw} to hold through every step, including the clipped step. We shall prove
\begin{equation}\label{eq:prefix}\textstyle E_W+p_W^2\le1\ \text{for every prefix $W$ of $F_n$.}\end{equation}
It suffices to check one insertion. Insert a list $Q$ satisfying the prefix bound between $P$ and $S$, where $P$ and $PS$ are prefixes of the original list. The new prefix $PQS$ has potential $E_P+p_P^2\{E_Q+p_Q^2(E_S+p_S^2)\}$. Using $E_Q\le1-p_Q^2$ bounds this by $(1-p_Q^2)(E_P+p_P^2)+p_Q^2(E_{PS}+p_{PS}^2)\le1$. A prefix ending inside $Q$ uses the same argument with $Q$ truncated and $S$ empty. Induction in \eqref{eq:table} transfers the raw bound to $F_n$. Finally the potential is convex in the last coefficient, and its values at zero and at the full coefficient are both at most one. This proves \eqref{eq:prefix}.

\subsection{A uniform bound on suffix products}\label{sec:suffix}
\begin{lemma}\label{lem:suffix}
Every suffix $W$ of $F_n$ satisfies $|p_W(z)|\le5$ for $0\le z\le1$.
\end{lemma}
\begin{proof}
Put $u=\arcsin\sqrt z$, $y=nu$ and $P=|p_{F_{n/4}}(z)|$ for $n\ge4$. We first identify which parts of a suffix can amplify its product. A suffix beginning in the last right copy reduces to a suffix at phase $y/4$. One beginning in the last left copy or at the last raw step contributes a multiplier at most $\max(1,|1-za_{n-1}|)P$ and reduces to the same phase. An earlier suffix contains both final copies, a raw suffix, and at most one incomplete initial copy of phase at most $y/16$. All other copies are complete and contract by \eqref{eq:prefix}. The multiplier in this case is at most $GP^2$, where $G=\max(1,y^2/2)$ bounds any raw suffix. Only its first partial dyadic layer can amplify it, and that layer has at most $n/2$ descendants of the root $2$, as in the proof of \Cref{lem:raw}. The starting degrees $n=1,2$ have only contractive factors. At $z=0$ all products equal one.

To bound the last-step multiplier, write $s(v)=\sin(v)/v$ with $s(0)=1$, and set
\[
\textstyle H(y)=(4y^2/\pi^2-1)\prod_{(d,k)\in\mathcal D}|s(y/d)|^k,\ \mathcal D=\{(2,1),(8,2),(32,6),(128,20)\}.
\]
We claim the multiplier is at most $\max(1,H(y))$. Writing $a=\csc^2 v$, the angles of its children are $(\pi-v)/2$ and $v/2$. Induction shows that each layer contains all odd multiples of its smallest angle below $\pi/2$. Thus $a_{n-1}=\csc^2(\pi/(2n))$ and every earlier angle is at least $\pi/n$. Put $v=\pi/(2n)$. If $u<v$ the last raw factor contracts. For $u\ge v$ we have
\[
\textstyle za_{n-1}-1\le (u^2/v^2-1)s(u)^2.
\]
Indeed $\csc^2 u-1/u^2$ increases on $(0,\pi/2]$ because $u^3\cos u\le\sin^3 u$. Indeed $u\cot u\le1-u^2/3$ because $(1-u^2/3)\sin u-u\cos u$ vanishes at zero and has derivative $u(\sin u-u\cos u)/3\ge0$. Also $s(u)\ge1-u^2/6$ gives $1-u^2/3\le s(u)^2$.

To bound $P$, recall from \Cref{lem:raw} that $p_{A_j}=(u\cot u)s(2ju)$. Inside $F_{n/4}$, the raw tables of degrees $n/4$, $n/16$, $n/64$ and $n/256$ occur with multiplicities $1$, $2$, $6$ and $20$, respectively, whenever the degree is at least one. Indeed the first three immediate child multiplicities at ratios $1/4,1/16,1/64$ are $2,2,4$ by \eqref{eq:table}, giving $6=2\cdot2+2$ and $20=2\cdot6+2\cdot2+4$. The other complete raw products have modulus at most one. This proves the claim for $n\ge256$. For the smaller degrees each retained raw table supplies a factor $u\cot u\le s(u)^2$, and the last raw step supplies another $s(u)^2$ when it amplifies. If $c\le1$ is a dyadic power, $s(u)\le s(cu)^{1/c^2}$ by iterating $s(2v)=s(v)\cos v\le s(v)^4$. Consequently $s(u)^b\le\prod_j s(c_j u)^{k_j}$ whenever $\sum_j k_j c_j^2\le b$. For $n\in\{4,8\},\{16,32\},\{64,128\}$ the available powers $b$ are $4,8,20$. The corresponding sums for the missing factors of $H$ are at most $2+6/16+20/256,\ 6+20/16,\ 20$, respectively, evaluated at the larger degree of each pair. This proves the claim for every $n$.

We next restrict the phases at which a multiplier can exceed one. If $y\le2$, monotonicity of $s$ on $[0,\pi/2]$ gives $za_{n-1}\le\max(1,4y^2/\pi^2)<2$. All raw factors therefore contract, and the same bound applies at each smaller inserted degree. For $y\le4$, every raw factor except the last contracts by the earlier-angle bound $za_j\le\max(1,y^2/\pi^2)<2$. Thus the earlier multiplier is also bounded by $\max(1,H(y))$ there. For $y>4$, $P\le2/y$ gives $GP^2\le2$. For $y\ge16$ necessarily $n\ge16$, and the first three complete raw tables in the guard give $P\le(2/y)(8/y)^2$, so $GP^2\le8192/y^4<1$. For $y\ge6\pi$ the last-step multiplier also contracts, since
\[
\textstyle H(y)\le \frac{512}{\pi^2 y}\sin^2(y/8)|s(y/32)|^6<1.
\]
On $[6\pi,8\pi]$ the last two factors are bounded by $1/2$ and $0.944^6$. On $[8\pi,32]$ their bounds are $3/4$ and $0.91^6$, using $|\sin(y/8)|\le4-\pi$. Beyond $32$ use $1$ and $(6/7)^6$. These bounds follow from monotonicity of $s$ up to $\pi$, $|s(v)|\le1/\pi$ thereafter, and $s(v)\le1-v^2/6+v^4/120$ at $v=3\pi/16,\pi/4,1$.

It suffices to bound the product of at most two amplifying multipliers. Indeed every reduction divides the phase by at least four, and three phases in $(2,6\pi)$ cannot lie on one path. Two such phases must be consecutive and equal to $y$ and $y/4$, with the first multiplier coming from the last raw step. For the second multiplier we can use $H(y/4)$ except in the earlier-suffix case with $y/4>4$, when the product is at most $2H(y)$ and $16<y<6\pi$. A single earlier multiplier can always be bounded by $2$ or $\max(1,H(y))$. The following five bounds finish all these cases. The entries in the last two columns are upper bounds.
\begin{center}
\begin{tabular}{c c r r c}
$J$ & interval & $t$ & $J(t)$ & $|(\log J)'(t)|$ \\ \hline
$H(y)$ & $[2,2\pi]$ & 4.132 & 2.27 & 0.001 \\
$H(y)$ & $[2\pi,4\pi]$ & 9.358 & 4.12 & 0.001 \\
$H(y)$ & $[4\pi,6\pi]$ & 15.017 & 2.25 & 0.001 \\
$H(y)H(y/4)$ & $[8,4\pi]$ & 10.206 & 4.39 & 0.001 \\
$H(y)H(y/4)$ & $[4\pi,6\pi]$ & 15.144 & 4.91 & 0.001
\end{tabular}
\end{center}
To check the entries by rational arithmetic, use $3.14159<\pi<3.14160$, obtained by taking the first seven terms of the alternating arctangent series in $\pi/4=4\arctan(1/5)-\arctan(1/239)$ and bounding each remainder by its next term. The logarithmic derivative is
\[
\textstyle (\log H)'(t)=2t/(t^2-\pi^2/4)+\sum_{(d,k)\in\mathcal D} k\bigl(\cot(t/d)/d-1/t\bigr).
\]
In each sine and cosine at $v=t/d$ or $t/(4d)$ substitute its Taylor polynomial at zero through degree $41$, with error interval $[-v^{42}/42!,v^{42}/42!]$. All arguments are below eight. For $J(y)=H(y)H(y/4)$ add $(\log H)'(t/4)/4$ to the displayed derivative. The entries bound whole intervals because $\log(y^2-\pi^2/4)$ is concave for $y>\pi/2$ and $(\log|s(v)|)''=v^{-2}-\csc^2v\le0$ between zeros. Thus at every interior point $J(y)\le J(t)\exp(0.001|y-t|)<1.008J(t)$, with endpoints following by continuity. The largest paired bound is $1.008\cdot4.91<5$, and the exceptional product is at most $2\cdot1.008\cdot2.25<5$. This proves the suffix bound.
\end{proof}

\subsection{Propagation of gradient errors}\label{sec:kernel}
\begin{lemma}\label{lem:kernel}
For every prefix of $F_n$, including one with its last coefficient shortened, the matrices in \eqref{eq:gains} satisfy
\[
\textstyle \sup_{0\le z\le1}\bigl(\|H(z)\|_F^2+\|H^\sharp(z)\|_F^2\bigr)\le (900)^2 n^4/2 .
\]
Consequently the nonlinear transfer estimate \eqref{eq:change} holds with $K=900$.
\end{lemma}
\begin{proof}
To control both gradient errors, consider any list $W=(h_0,\ldots,h_{t-1})$ and fix $z\in[0,1]$. Write $p_{<i},p_{>i},p_{\ge i},p_{(s,i)}$ for the products of $1-zh_j$ over the indicated indices. Directly from \eqref{eq:gains},
\[
\textstyle G(W):=\|H\|_F^2+\|H^\sharp\|_F^2=\sum_{s<i}h_s h_i(1+zh_i)p_{(s,i)}^2.
\]
It suffices to bound $G(F_n)$, since truncation deletes summands and shortening the last coefficient can only reduce the last row. Partition $F_n$ into its raw singleton steps and its nonempty inserted copies. The terms with both indices in one block contribute $2\sum_j G(F_{r_j})$, where $G(F_0)=0$. For distinct blocks $X$ before $Y$, the contribution factors as $T(X)p_I^2 U(Y)$, where $I$ consists of the complete intervening blocks and
\[
\textstyle T(W)=\sum_i h_i p_{>i}^2,\ U(W)=\sum_i h_i(1+zh_i)p_{<i}^2.
\]
We first bound the endpoint sums and then the intervening product.

The mass satisfies $M(F_n)\le n^2$. Indeed $M(F_1)=1$, and for $n\ge2$ the raw mass is at most $3n^2/4$ and the inserted mass at most $2\sum_jr_j^2\le15n^2/112$ by induction. Put $x=zn^2$. When $x\le1$, all factors have modulus at most one and $zh_i\le1$, giving $T\le n^2$, $U\le2n^2$, and $G(F_n)\le2\sum_{s<i}h_sh_i\le n^4$. For $x>1$, telescoping the forward products gives $z^2\sum_i h_i^2p_{<i}^2=2E_{F_n}+p_{F_n}^2-1\le E_{F_n}$ by \eqref{eq:prefix}. Hence $U(F_n)\le2\min(n^2,1/z)$ for $z>0$.

For the other endpoint sum we need $T(F_n)\le 284 n/\sqrt z$ and $T(F_n)\le18952\min(n^2,1/z)$. To obtain both bounds, set $y_i=zh_i$ and telescope the suffix products in reverse order, obtaining
\[
\textstyle zT(F_n)=1-p_{F_n}^2+\sum_i y_i(y_i-1)p_{>i}^2.
\]
Only indices with $y_i>1$ can increase this bound. Their contribution is at most $2\max(p_{>i}^2,p_{\ge i}^2)$ per index, since $y(y-1)\le2$ for $1<y\le2$ and $y(y-1)\le2(y-1)^2$ for $y\ge2$. There are at most $4\sqrt{2}\,n\sqrt z$ such indices. Indeed the raw coefficient bound $a_j\le2q_j^2$ requires $q_j>1/\sqrt{2z}$, so $j+1$ is a multiple of the smallest dyadic scale above $1/\sqrt{2z}$, giving at most $\sqrt{2}\,n\sqrt z$ raw indices. Induction in the inserted copies gives at most $4\sqrt{2z}\,(2\sum_j r_j)\le3\sqrt{2z}\,n$ more. The bound $5$ of \Cref{lem:suffix} therefore gives $zT\le1+8\sqrt{2}\cdot 25\sqrt x$. For $x\ge1$ this proves $T\le 284 n/\sqrt z$. For $x\le1$ use $T\le n^2$.

For a bound on $zT$ independent of $x$, group the same indices along the chain of final copies $F_n,F_{n/4},\ldots$. At degree $d\ge4$ take the indices outside its final copy and descend into that copy. At degree one or two take all remaining indices. Along this chain the suffixes of any retained index agree with its suffixes in $F_n$. Write $u=zd^2$. Both suffixes of an index outside the final copy contain that copy. Each contains at most one incomplete initial insertion, complete insertions, and a raw suffix. Any raw suffix has modulus at most $\max(1,u/2)$, since it meets a suffix of a layer with at most $d/2$ descendants of the coefficient $2$, followed by complete layers of modulus at most one. This uses the descendant interval bound in the proof of \Cref{lem:raw}. The initial unit, if present, is contractive. If $u>2048$, then $d\ge64$. The final copy contains the raw products of one $A_{d/4}$, two $A_{d/16}$ and at least four $A_{d/64}$, by applying the two insertions at the final raw site twice. Each additional complete raw product has modulus at most one. Thus \Cref{lem:raw} bounds the two suffix moduli by $5\eta(u)$, where
\[
\textstyle \eta(u)=1\ (u\le2048),\ \eta(u)=(64/\sqrt{u})(32/\sqrt{u})^4\ (u>2048).
\]
Here the large-$u$ bound multiplies $u/2$ by $(2/\sqrt{u})(8/\sqrt{u})^2(32/\sqrt{u})^4$. For $u\le2048$ we use the uniform suffix bound directly. The scales $u$ decrease by a factor of $16$. Thus the sum of $\sqrt{u}\,\eta(u)^2$ is less than $67$. The terms with $u\le2048$ sum to at most $(4/3)\sqrt{2048}<61$, and for $u>2048$ the summand is $4096(1024/u)^4/\sqrt u$, giving a sum at most $(256/\sqrt{2048})/(1-16^{-9/2})<6$. Using the index count at each degree in the reverse telescoping identity yields $zT\le1+8\sqrt{2}\cdot 25\cdot67<18952$. Together with $T\le n^2$ for $x\le1$, this proves the second bound.

We now control the terms from distinct blocks when $x\ge1$. Give a raw step $a$ weight $w=\min(a,1/z)$ and a copy $F_r$ weight $w=\min(r^2,1/z)$. Let $W_R,W_I$ be the two total weights. The scale $q$ occurs $n/(2q)$ times with an additional $1/2$ at $q=n$. Using $a\le2q^2$ and splitting at $q=1/\sqrt z$, the scales below and above the split and the extra term contribute at most $2n/\sqrt z$, $n/\sqrt z$ and $n/(2\sqrt z)$, respectively. Thus $W_R\le4n/\sqrt z$. For the two copies at each site use $r\le q/4$ and split $\sum_q(n/q)\min(q^2/16,1/z)$ at $4/\sqrt z$. Each side sums to at most $n/(2\sqrt z)$, and the extra contribution $\min(n^2/16,1/z)$ is at most $n/(4\sqrt z)$. Hence $W_I\le5n/(4\sqrt z)$.

To bound $p_I$ together with the endpoint factors, include raw endpoints when locating the largest raw scale $q$ between the two blocks. Its site is unique, since a larger dyadic scale lies between two equal ones. A left raw endpoint has $\sqrt{a/w}\le 1+|1-za|/2$, and a right raw endpoint has $\sqrt{a(1+za)/w}\le 3/2+|1-za|$. Expanding the products of these bounds either includes or omits each raw endpoint. Each resulting raw interval has product of modulus at most $\max(1,2zq^2)$ by the proof of \Cref{lem:raw}. If a copy at the largest site is strictly between the endpoints, its product has modulus at most $\min(1,1/(2r\sqrt z))$ by \Cref{lem:raw} and \eqref{eq:prefix}. Their product is at most $4\sqrt x$ because $q^2=4nr$ and $x\ge1$. When that site has no copy, $q^2\le2n$ by dyadicity, so the raw bound alone is at most $4\sqrt x$. A pair containing no raw step between or at its endpoints also obeys this bound since all intervening copies contract. The sums of the expansion coefficients for a left and a right raw endpoint are $3/2$ and $5/2$. Squaring, summing over these pairs, and using $2\sum_j r_j\le3n/4$ now bounds their total by
\[
\textstyle 16x\bigl[(9/4)W_R+\sum_{\rm copies}T(F_r)\bigr]\bigl[(25/4)W_R+2W_I\bigr]\le16(9+3\cdot284/4)(25+5/2)n^4 .
\]

For the remaining pairs neither flanking copy of the largest site lies strictly between the endpoints. Both endpoints therefore belong to the same triple $(F_r,(a),F_r)$. Write $w=\min(q^2,1/z)$ and $C_T=18952$. Since $a\le2q^2$ and $T(F_r)\le C_Tw$, $U(F_r)\le2w$, the left copy and the raw step contribute at most $6C_Tq^4$, using $(1+2zq^2)w\le3q^2$. The raw step and the right copy contribute at most $4q^4$. For the two copies the intervening raw factor gives at most $8C_Tq^4$, since $w\max(1,2zq^2)\le2q^2$. The sum of $q^4$ over sites is at most $(8/7)n^4$ by the geometric series with at most $n/q$ sites of scale $q$. Finally $r\le q^2/(4n)$ at each site gives $2\sum_j r_j^4\le (2n^4/4^4)\sum_{k\ge0}2^{-7k}=n^4/127$. Including the within-block pairs, induction on $n$ gives
\[
\textstyle G(F_n)/n^4\le\tfrac{127}{126}\bigl[16(9+3\cdot284/4)(25+5/2)+(8/7)(14\cdot18952+4)\bigr] < (900)^2/2.
\]
For the final assertion, apply $H(z)$ and $H^\sharp(z)$ to the same input on each positive eigendirection of the frozen Hessian. The sum of their squared error norms is bounded by $(\sup_z G)$ times the total squared input norm. Cauchy--Schwarz bounds the sum of the two error norms by $\sqrt{2\sup_z G}$ times that input norm. This gives the coefficient $900n^2$ used in \eqref{eq:change}.
\end{proof}

\section{Sufficient Decrease Within an Epoch}\label{app:descent}
For the proofs in this and the following appendices, write $A=\sum_{j=0}^{J-1}b_j$, $S_t=\sum_{j<t}h_j\norm{g_j}^2$ and $C_t=\sum_{j<t}h_j\inner{g_j}{g_{j+1}}$. The test in Algorithm~\ref{alg:gd} halves $B$ exactly when $S_t+2C_t<0$. At a completed epoch the amount spent equals $T=\max\{B,\norm{g_0}^2\}$, since only the first, unit step can spend more than $B$.

\paragraph{Two estimates for the epoch proof.} For a nonempty executed prefix $x_0,\ldots,x_t$ with cap $T=B$ let $R=\sqrt{A S_t}/L$. Cauchy--Schwarz gives $\sum_{i<t}\norm{x_{i+1}-x_i}\le R$. We need one estimate to translate the observable test into a decrease and another to show that the test eventually passes. We state them together so that their roles in \Cref{prop:epoch} are explicit.
\begin{lemma}\label{lem:descent}
For each nonempty prefix with $T=B$, the following bounds hold with $K=900$.
\begin{equation}\label{eq:trap}\textstyle \bigl|\tfrac12(S_t+C_t)-L(f(x_0)-f(x_t))\bigr|\le\tfrac{\rho R}{12L}\max_{i<t}h_i\, S_t,\end{equation}
\begin{equation}\label{eq:change}\textstyle f(x_t)-f(x_0)\le\tfrac{S_t}{L}\bigl[-1+\tfrac12(1+K m^2\rho R/L)^2+\tfrac{\rho R}{2L}\max_{i<t}h_i\bigr].\end{equation}
\end{lemma}

\begin{proof}[Proof of \Cref{prop:epoch}]
To bound the decrease in a passing epoch we must lower-bound its actual decrease from $C_t\ge-S_t/2$. First dispose of $T>B$. Here $h_0=1$ and $\norm{g_1-g_0}\le\norm{g_0}$ gives $\inner{g_0}{g_1}\ge0$. The test passes and the fixed-step bound gives decrease at least $S_1/(2L)\ge B/(2L)$. For $T=B$ we have $S_t=B$ at completion. The table bounds $A\le2m^2$ and $\max h_i\le m^2/2$ from Appendix~\ref{sec:stability}, together with $B_0=\eps^2A/2$ and $\rho\eps m^4/L^2\le1$, give
\begin{equation}\label{eq:scale}\textstyle \tfrac{\rho R}{L}\max h_i\le\tfrac1{\sqrt2}\sqrt{B/B_0},\ m^2\tfrac{\rho R}{L}\le\sqrt2\,\sqrt{B/B_0}.\end{equation}
Consequently \eqref{eq:trap} bounds the error by $B/(12\sqrt2)<B/16$. If the test passes, its estimate $(S_t+C_t)/2$ is at least $B/4$, proving the claimed decrease of $3B/(16L)$.

It remains to control a failure and show that it requires $B$ to exceed a fixed fraction of $B_0$. Since $B/B_0\le1$, \eqref{eq:scale} bounds the bracket in \eqref{eq:change} by an absolute constant. The increase in $f$ is therefore $\cO(B/L)$. As $B/B_0$ decreases, the same bracket is bounded above by an expression tending to $-1/2$, uniformly over epochs. Choose an absolute $c>0$ small enough that it is at most $-3/8$ for $B/B_0\le c$. The actual decrease is then at least $3B/(8L)$, and \eqref{eq:trap} yields $(S_t+C_t)/2\ge 3B/8-B/16>B/4$. The test must pass.
\end{proof}

\begin{proof}[Proof of \Cref{lem:descent}]
For \eqref{eq:trap}, our goal is to bound the discrepancy between the gradient test and the true decrease. On each segment the decrease is the integral of the negative directional derivative, while the test uses its two endpoint values. Writing $s_i=x_{i+1}-x_i$, the derivative of $u\mapsto\inner{\nabla f(x_i+us_i)}{s_i}$ is Lipschitz with constant $\rho\norm{s_i}^3$. The unit-interval trapezoidal error is at most this constant divided by $12$, by integrating the kernel $u(1-u)/2$. Multiplication by $L$ and summation give an error at most $(L\rho/12)\sum_{i<t}\norm{s_i}^3$. Now $\max_i\norm{s_i}\le R$ and $\sum_i\norm{s_i}^2\le S_t\max_i h_i/L^2$, which proves \eqref{eq:trap}.

For \eqref{eq:change}, the linear terms of Taylor's inequality already sum to $-S_t/L$. It remains to control its quadratic terms. Freeze $M=\nabla^2 f(x_0)$ and let $(M/L)_+$ denote the matrix obtained by replacing its negative eigenvalues by zero. Every Hessian on the segments is bounded above by $M+\rho R I$. Taylor's inequality thus gives
\begin{equation}\label{eq:taylor}\textstyle f(x_t)-f(x_0)\le-\tfrac{S_t}{L}+\tfrac1{2L}\sum_{i<t} h_i^2\inner{g_i}{(M/L)_+g_i}+\tfrac{S_t}{2L}\tfrac{\rho R}{L}\max_{i<t}h_i.\end{equation}
We must bound the middle sum by $(1+Km^2\rho R/L)^2S_t$. We do this first for gradients predicted by the frozen Hessian and then bound the error of that prediction.

On an eigendirection of $M/L$ of eigenvalue $z\in(0,1]$, the frozen prediction at step $i$ is $p_i(z)g_0(z)$, where $p_i(z)=\prod_{j<i}(1-h_jz)$ and $g_i(z)$ is the coordinate in that direction. The scalar prefix estimate \eqref{eq:prefix} gives $z\sum_{i<t}h_ip_i^2+p_t^2\le1$. Telescoping $p_{i+1}^2-p_i^2=(-2zh_i+z^2h_i^2)p_i^2$ therefore gives
\begin{equation}\label{eq:frozen}\textstyle z\sum_{i<t}h_i^2p_i(z)^2\le\sum_{i<t}h_ip_i(z)^2.\end{equation}
Thus for the predictions the desired weighted norm is at most the ordinary norm with weights $h_i$. To transfer this comparison to the actual gradients, we need to bound the prediction error in both norms.

The error injected in step $s$ is $\sqrt{h_s}w_s=g_{s+1}-g_s+(h_s/L)Mg_s$. Integrating the Hessian along that step gives $\norm{w_s}\le(\rho R/L)\sqrt{h_s}\norm{g_s}$. Thus $\sum_{s<t}\norm{w_s}^2\le(\rho R/L)^2 S_t$. Iterating the gradient recurrence, the vector with entries $\sqrt{h_i}(g_i(z)-p_i(z)g_0(z))$ is $H(z)w(z)$, where $w(z)=(w_s(z))_{s<t}$. Multiplication of entry $i$ by $\sqrt{zh_i}$ gives $H^\sharp(z)w(z)$. The matrices are zero for $i\le s$ and otherwise have entries
\begin{equation}\label{eq:gains}\textstyle H_{is}(z)=\sqrt{h_i h_s}\prod_{s<j<i}(1-zh_j),\ H^\sharp_{is}(z)=\sqrt{zh_i}\,H_{is}(z).\end{equation}
By \Cref{lem:kernel}, summing over an orthonormal basis of the positive eigenspace bounds the sum of these two error norms by $900m^2(\rho R/L)\sqrt{S_t}$, also when the final step is shortened. The ordinary weighted norm of the actual gradients in this eigenspace is at most $\sqrt{S_t}$. Use the triangle inequality first to pass from the actual weighted curvature norm to its prediction, then use \eqref{eq:frozen}, and use the triangle inequality again to return to the ordinary weighted norm of the actual gradients. This gives
\[\textstyle \bigl(\sum_{i<t}h_i^2\inner{g_i}{(M/L)_+g_i}\bigr)^{1/2}\le(1+Km^2\rho R/L)\sqrt{S_t}.\]
Substituting in \eqref{eq:taylor} proves \eqref{eq:change}. All spectral decompositions here use the single symmetric matrix $M$.
\end{proof}

\section{Gradient Evaluation Complexity}\label{app:complexity}
\begin{proof}[Proof of \Cref{thm:main}]
To obtain \eqref{eq:target} we need the cost of one epoch and the number of epochs. For the cost, Appendix~\ref{sec:stability} proves $J\le4m$ and $(2m^2+1)/3\le A\le2m^2$. If $T>B$, the unit first step spends the entire cap. If $T=B$, spending the full table while all gradients exceed $\eps$ would cost more than $\eps^2A=2B_0\ge2B$. Thus each epoch reaches the cap or returns an iterate within $J$ updates.

For the count, \Cref{prop:epoch} implies that each failure halves $B$ and that failures stop once $B\le cB_0$. Hence there are at most $\lceil\log_2(1/c)\rceil$ failures, $B$ stays at least $cB_0/2$, and its values just before the halvings sum to at most $2B_0$. All possible increases therefore total $\cO(B_0/L)$. Each passing epoch decreases $f$ by at least $3cB_0/(32L)$. Telescoping through the last completed nonterminal epoch and using $f\ge\inf f$ bounds the number of such epochs by $\cO(1+L\Delta/B_0)$. Including the last epoch and the initial gradient gives $\cO(1+m+L\Delta/(\eps^2 m))$ evaluations. When $m\ge2$, its definition gives $m\asymp(L^2/(\rho\eps))^{1/4}$ and proves \eqref{eq:rate}. When $m=1$, $B_0=\eps^2/2$ and every nonterminal step is ordinary $1/L$ gradient descent. The fixed-step bound from \eqref{eq:fixed} suffices because $L\Delta/\eps^2=(L^2/(\rho\eps))^{1/4}\Delta L^{1/2}\rho^{1/4}\eps^{-7/4}$ and $(L^2/(\rho\eps))^{1/4}<2$.
\end{proof}

\section{Explicit gradient evaluation bounds}\label{app:explicit}
\subsection{Explicit descent at degrees four and eight}\label{sec:small-degrees}
We show that every completed epoch passes at its initial threshold $B_0$ when Algorithm~\ref{alg:gd} selects $m=4$ or $m=8$. We then count the steps, including a possible last epoch that stops early. Write $a=A/2$ and let $W=(h_0,\ldots,h_{J_m-1})$ be the nominal table cut off after coefficient sum $a$, shortening its last coefficient to attain $a$. Write $h_{\max}=\max W$. If $T>B_0$ the unit step passes and decreases $f$ by at least $S_1/(2L)$. Otherwise $T=B_0=a\eps^2$. Each source gradient exceeds $\eps$, so $\sum_{i<t}h_i\le S_t/\eps^2\le a$ for the actual steps. They form a prefix of $W$ with at most its final coefficient shortened, and
\begin{equation}\label{eq:shortpath}\textstyle \sum_{i<t}\|x_{i+1}-x_i\|=L^{-1}\sum_{i<t}h_i\|g_i\|\le S_t/(L\eps).\end{equation}

To bound the actual decrease, freeze the Hessian at the midpoint of this polygonal path measured in arclength, and call it $M$. Along the path the normalized Hessian differs from $M/L$ by at most $\delta=\rho S_t/(2L^2\eps)\le a/(2m^4)$. Taylor's inequality gives
\[
\textstyle L(f(x_0)-f(x_t))\ge S_t-\tfrac12\sum_{i<t}h_i^2\inner{g_i}{(M/L)_+g_i}-\tfrac12\delta h_{\max}S_t.
\]
To bound the middle sum, put $w_i=(g_{i+1}-g_i+(h_i/L)Mg_i)/\sqrt{h_i}$. Integrating the Hessian on each step gives $\sum_{i<t}\|w_i\|^2\le\delta^2 S_t$. At a positive eigenvalue $z$ of $M/L$, the corresponding scalar coordinates satisfy $v_{i+1}=(1-zh_i)v_i+\sqrt{h_i}\,w_i(z)$. It suffices to prove
\begin{equation}\label{eq:storage-goal}\textstyle \sum_{i<t}(zh_i^2-\beta h_i)v_i^2\le m^4\sum_{i<t}|w_i(z)|^2,\ \beta=21/20.\end{equation}
Indeed, summing over positive eigendirections bounds the middle sum by $(\beta+m^4\delta^2)S_t$. We then get $f(x_0)-f(x_t)\ge d_m S_t/L$ with
\begin{equation}\label{eq:delta-exp}\textstyle d_m=\tfrac{19}{40}-\tfrac{a^2}{8m^4}-\tfrac{a h_{\max}}{4m^4}.\end{equation}
Using \eqref{eq:shortpath} in \eqref{eq:trap} bounds the trapezoidal error by $S_t a h_{\max}/(12m^4)$. Thus $d_m-a h_{\max}/(12m^4)\ge1/4$ implies passage.

To establish \eqref{eq:storage-goal} simultaneously for all prefixes, it suffices to find nonnegative $R_0,\ldots,R_{J_m}$ with $R_0=R_{J_m}=0$ such that at each nominal coefficient $h_i$ of $W$,
\begin{equation}\label{eq:ric}\textstyle h_iR_{i+1}<m^4,\ R_i\ge zh_i^2-\beta h_i+\frac{(1-zh_i)^2R_{i+1}}{1-h_i R_{i+1}/m^4}.\end{equation}
Completing the square in $w$ bounds $R_{i+1}((1-zh_i)v+\sqrt{h_i}w)^2-m^4w^2$ by the last fraction times $v^2$. Summing the resulting inequalities telescopes. If the last coefficient is shortened to $s\le h_i$, convexity gives $zs^2-\beta s\le\max(0,zh_i^2-\beta h_i)\le R_i$, so the last row still gives \eqref{eq:storage-goal}.

For an interval $[l,r]$ of curvatures and bounds $u\le h_i\le v$, the right side of \eqref{eq:ric} is at most
\begin{equation}\label{eq:finite-row}\textstyle \max_{s\in\{u,v\}}(rs^2-\beta s)+\frac{R_{i+1}\max\{(1-lu)^2,(1-rv)^2\}}{1-vR_{i+1}/m^4}.\end{equation}
The first term uses convexity in $s$, and $1-zh_i$ lies between $1-rv$ and $1-lu$. For $m=4$, the truncated table is $(1,2,h,1,5/2-h)$, where $h=4-2\sqrt2\in[7/6,6/5]$ and $5/2-h\in[13/10,4/3]$. With these bounds, substitution in \eqref{eq:finite-row} gives the following choices.
\begin{center}
\begin{tabular}{c|rrrr}
$[l,r]$ & $R_1$ & $R_2$ & $R_3$ & $R_4$\\ \hline
$[0,1/2]$ &0&0&0&0\\
$[1/2,3/4]$ &1&0&0&0\\
$[3/4,7/8]$ &$7/5$&0&0&$1/5$\\
$[7/8,1]$ &$9/4$&$1/5$&0&$2/5$
\end{tabular}
\end{center}

For $m=8$, the truncated table is $(1,2,a_2,a_3,a_4,a_5,a_6,1,2,21/2-4a_2-a_6)$. Here $a_3=8-a_2$ and $a_5=4a_2-a_4$. To bound the radical terms using rationals, the smaller child $s$ of a parent $p$ satisfies $s^2/(4(s-1))=p$, with the left side decreasing on $(1,2)$. Substitution of endpoints, first for $p=2$, then for $p=a_2,a_3$, gives $1.17157\le a_2\le1.17158$, $1.4464\le a_4\le1.4465$, and $1.0395\le a_6\le1.0396$. The resulting lower and upper bounds for $100h_i$, in order, are
\begin{center}
\begin{tabular}{c|rrrrrrrrrr}
$i$&0&1&2&3&4&5&6&7&8&9\\ \hline
lower&100&200&117&682&144&323&103&100&200&477\\
upper&100&200&118&683&145&324&104&100&200&478
\end{tabular}
\end{center}
The following table gives the nine interior $R_i$ using these bounds. In every row the denominators in \eqref{eq:finite-row} are positive, since $6.83\cdot142<8^4$.
\begin{center}
\begin{tabular}{c|rrrrrrrrr}
$[l,r]$ & $R_1$ & $R_2$ & $R_3$ & $R_4$ & $R_5$ & $R_6$ & $R_7$ & $R_8$ & $R_9$\\ \hline
$[0,1/8]$ &0&0&0&0&0&0&0&0&0\\
$[1/8,1/4]$ &1&3&5&0&0&0&0&0&1\\
$[1/4,3/8]$ &1&5&11&0&1&0&0&1&4\\
$[3/8,1/2]$ &1&5&17&0&2&0&0&1&7\\
$[1/2,5/8]$ &1&6&33&1&4&0&1&2&10\\
$[5/8,3/4]$ &2&4&45&1&5&0&1&5&13\\
$[3/4,7/8]$ &2&1&59&1&6&0&1&10&15\\
$[7/8,15/16]$ &4&2&96&2&7&0&1&15&17\\
$[15/16,1]$ &7&5&142&3&8&0&1&21&18
\end{tabular}
\end{center}
For example, the last row at $i=1$ bounds \eqref{eq:finite-row} by $4-2\beta+5/(1-10/4096)<7$. At $i=0$ the upper bound is $1-\beta+(7/256)/(1-7/4096)<0$, as required by $R_0=0$. Every interior entry is obtained by taking the nonnegative part of \eqref{eq:finite-row} and rounding up to an integer, starting at $R_{10}=0$.

It remains to count the evaluations. The pairs $(a,J_m)$ are $(13/2,5)$ and $(49/2,10)$. Substituting $h_{\max}\le2$ and $h_{\max}\le6.83$, respectively, in \eqref{eq:delta-exp} gives $d_4\ge18091/40960$ and $d_8\ge365747/819200$, with $d_m-a h_{\max}/(12m^4)>0.43$ in both cases. Thus each completed epoch passes, and induction keeps $B=B_0$ throughout. A capped epoch takes at most $J_m$ updates and decreases $f$ by at least $d_ma\eps^2/L$. Its number of updates divided by its decrease in units of $\eps^2/L$ is at most $J_m/(a d_m)$, which is below $1.742$ for $m=4$ and $0.915$ for $m=8$.

Only the final epoch can stop strictly before its cap. If it takes $t$ updates, their coefficients are nominal and $t<J_m$. Its decrease in the same units is at least $d_m\sum_{i<t}h_i$, or at least $1/2$ for $t=1$. For $m=4$ the first four prefix sums are bounded below by $(100,300,417,517)/100$, and for $m=8$ the first nine are bounded below by $(100,300,417,1099,1243,1566,1669,1769,1969)/100$. Multiplying these bounds by $d_m$, using $1/2$ at $t=1$, and subtracting $1.742$ or $0.915$ times the resulting decrease from $t$ leaves at most $0.129$ for $m=4$ and $1.30$ for $m=8$. Adding the capped epochs and the initial gradient evaluation gives $N\le1+1.742X+0.129$ for $m=4$ and $N\le1+0.915X+1.30$ for $m=8$, where $X=L\Delta/\eps^2$.

The same prefix sums show $\sum_{i<t}h_i\ge5t/4$ for $2\le t<J_m$. Since both $d_m>11/25$ and both capped charges are below $20/11$, each capped epoch and each strict terminal prefix with at least two steps costs at most $20/11$ times its decrease in these units. A one-step terminal epoch costs at most twice its decrease and needs at most $1/11$ extra at the smaller multiplier. Thus $N\le1+\min\{2X,(20X+1)/11\}$ at both degrees.

\subsection{Explicit bounds for arbitrary degree}\label{sec:all-degrees}
Write $X=L\Delta/\eps^2$ and $r=(L^2/(\rho\eps))^{1/4}$. We first quantify the guarantee for the degree in Algorithm~\ref{alg:gd}. We then give a choice of degree that ensures every completed epoch passes at its initial threshold.

\begin{proposition}\label{prop:original-count}
For the degree $m\ge2$ selected by Algorithm~\ref{alg:gd}, the number of gradient evaluations satisfies $N\le1+1.13\cdot10^7m+770X/m$.
\end{proposition}
\begin{proof}
We bound the number of updates in passing epochs by their decrease, then account for the increases and the updates in failed epochs. Put $e=B/\eps^2$. If $T=B$, the path length up to each prefix is at most $S_t/(L\eps)$ by \eqref{eq:shortpath}. Freeze the Hessian at the arclength midpoint of this path. Its normalized difference from every segment Hessian is bounded by $\delta=\rho S_t/(2L^2\eps)\le e/(2m^4)$, since $m\le r$. The transfer argument in \eqref{eq:taylor}--\eqref{eq:gains} uses only this difference, so \Cref{lem:kernel} gives
\begin{equation}\label{eq:both-midpoint}\textstyle
\tfrac{L(f(x_0)-f(x_t))}{S_t}\ge1-\tfrac12\bigl[(1+900m^2\delta)^2+h_{\max}\delta\bigr],\
\tfrac{C_t}{S_t}\ge1-(1+900m^2\delta)^2-h_{\max}\delta .
\end{equation}
For the second inequality, let $\bar M_i$ be the Hessian averaged along step $i$. Then $C_t=S_t-\sum_{i<t}h_i^2\langle g_i,(\bar M_i/L)g_i\rangle$, and the same curvature estimate applies. Here $h_{\max}\le m^2/2$. If $e\le m^2/2250$, the sum subtracted in the second inequality is at most $(1+1/5)^2+1/9000<3/2$. Every such completed epoch passes. A unit epoch with $T>B$ always passes by \Cref{prop:epoch}. Since $m^2/3\le B_0/\eps^2\le m^2/2$, it follows that $e\ge m^2/4500$ and there are at most $\lceil\log_2 1125\rceil=11$ failures.

For the cost of a passing epoch with $e\ge1$, let $k$ be the smallest dyadic degree whose raw mass $(2k^2+1)/3$ is at least $e$. It exists with $k\le m$, and $k\le\sqrt{6e}$ by minimality. The gradients at executed steps exceed $\eps$ in norm, so the threshold is reached no later than raw index $k-1$. Counting all insertions through that site except its right copy gives at most $k+8\sum_{j<k-1}r_j+4r_{k-1}$ steps, since $N(F_d)\le4d$. Each earlier insertion satisfies $r_j\le q_j^2/(4m)$, the last satisfies $r_{k-1}\le k^2/(4m)$, and $\sum_{j<k}q_j^2=(3k^2-k)/2$ by the dyadic site count. The step count is thus at most $k+2k^2/m\le\sqrt{6e}+12e/m$. The trapezoidal error is at most $S_t e h_{\max}/(12m^4)\le S_t/48$ by \eqref{eq:shortpath}, because $e\le m^2/2$. Hence a passing epoch decreases $f$ by at least $11B/(48L)$, also for $T>B$ by ordinary unit-step descent. Its update count is therefore at most $770/m$ times its decrease in units of $\eps^2/L$, since $(48/11)(\sqrt{27000}+12)<770$. If $e<1$, a single unit step completes the epoch and decreases $f$ by more than $\eps^2/(2L)$. Here $m<\sqrt{4500}$, so the same charge applies.

It remains to bound the total increase in failed epochs. By the first inequality in \eqref{eq:both-midpoint}, the increase in units of $\eps^2/L$ is at most $(450+1/8)e^2/m^2+(900^2/8)e^3/m^4$, after dropping the negative term $-e/2$. The successive values of $e$ at failures halve, starting at most at $m^2/2$. Summing the two geometric series bounds all increases by $[(450+1/8)/3+900^2/56]m^2<14615m^2$. Telescoping through the last completed nonterminal epoch bounds the total passing decrease in the same units by $X+14615m^2$. The eleven possible failed epochs and the terminal epoch require at most $48m$ further updates by $N(F_m)\le4m$. Including the initial evaluation gives the claimed count.
\end{proof}

\begin{corollary}\label{cor:explicit-rate}
Suppose $r=(L^2/(\rho\eps))^{1/4}\ge16$. Run Algorithm~\ref{alg:gd} with $m$ the largest power of two at most $r/8$. Every completed epoch passes at $B=B_0$, and
\begin{equation}\label{eq:explicit-rate}\textstyle
N\le1+3m+21X/m\le1+\tfrac38(L^2/(\rho\eps))^{1/4}+336\Delta L^{1/2}\rho^{1/4}\eps^{-7/4}.
\end{equation}
\end{corollary}
\begin{proof}
To count epochs without any failures, it suffices to prove descent and passage at $B_0$. Here $m\ge2$, $\rho\eps m^4/L^2\le8^{-4}$ and $B_0/\eps^2\le m^2/2$. At the path midpoint this gives $m^2\delta\le1/16384$ and $h_{\max}\delta\le1/32768$. Since $(1+900/16384)^2+1/32768<28/25$, \eqref{eq:both-midpoint} gives a decrease of at least $11S_t/(25L)$ and $C_t\ge-3S_t/25$. Unit epochs with $T>B_0$ satisfy both the decrease and the passage test. Induction now keeps $B=B_0$ throughout the run.

Each epoch takes at most $3m$ updates. Indeed, for $m\ge4$ the copy after the last raw step has mass at most $(m/4)^2<A/2$. Thus the mass through that step exceeds $A/2$, so the epoch has already stopped or reached $B_0$. Counting the steps through it gives $m+4(2\sum_jr_j-m/4)\le3m$. For $m=2$ the table has only two steps. Each completed nonterminal epoch decreases $f$ by at least $11B_0/(25L)$, where $B_0\ge m^2\eps^2/3$. Their number is therefore at most $75X/(11m^2)$. Allowing $3m$ updates for the terminal epoch and one initial evaluation proves the first bound, since $225/11<21$. The second uses $r/16<m\le r/8$.
\end{proof}

For the choice in \Cref{cor:explicit-rate} we can also compare with the fixed-step bound when $m\ge16$. Since $21/m<20/11$, the preceding proof already charges every completed nonterminal epoch, and every terminal epoch reaching its cap, at most $20/11$ times its decrease in units of $\eps^2/L$. To bound a strict terminal epoch with $t\ge2$, it suffices to show that its nominal prefix has mass at least $5t/4$, for then its decrease in the same units is at least $(11/25)(5t/4)=11t/20$.

Subtract $5/4$ from each coefficient and call the resulting sum its surplus. Induct on the table degree. A prefix of a smaller inserted copy has surplus at least $-1/4$, since by induction only its first step can have negative surplus. The first raw pair $(1,2)$ has no insertions and yields surplus $1/2$. Each later raw pair has mass $4a_j\ge4$, and only its second step can have insertions. A complete pair with its two copies therefore adds surplus at least $4-5/2-2/4=1$. Until the second raw step is taken, the first step and a prefix of the first copy each cost at most $1/4$, covered by the preceding surplus. Once the second step is taken, the same full-pair bound applies even with a partial final copy. This proves the prefix bound, starting with $F_1$ and $F_2$. A single terminal unit step decreases $f$ by at least $\eps^2/(2L)$. As for degrees four and eight above, summation yields $N\le1+\min\{2X,(20X+1)/11\}$.
\end{appendices}

\end{document}

%% file: header.tex
\usepackage[margin=1in]{geometry}
\usepackage{amsmath,amssymb,amsthm,mathtools,bm}
\usepackage{xcolor}
\definecolor{EvenRed}{HTML}{8D263A}
\definecolor{OddBlue}{HTML}{3D647D}
\definecolor{PlotGray}{HTML}{7A858C}
\usepackage{algorithm}
\usepackage{algpseudocode}
\usepackage{placeins}
\usepackage[linktocpage,colorlinks,linkcolor=blue,anchorcolor=blue,citecolor=blue,urlcolor=blue]{hyperref}
\usepackage[authoryear,round]{natbib}
\usepackage[title,titletoc]{appendix}

\usepackage{aliascnt}

\theoremstyle{plain}
\newtheorem{theorem}{Theorem}[section]

\newaliascnt{lemma}{theorem}
\newtheorem{lemma}[lemma]{Lemma}
\aliascntresetthe{lemma}

\newaliascnt{proposition}{theorem}
\newtheorem{proposition}[proposition]{Proposition}
\aliascntresetthe{proposition}

\newaliascnt{corollary}{theorem}
\newtheorem{corollary}[corollary]{Corollary}
\aliascntresetthe{corollary}

\theoremstyle{definition}

\newaliascnt{assumption}{theorem}
\newtheorem{assumption}[assumption]{Assumption}
\aliascntresetthe{assumption}

\newaliascnt{definition}{theorem}
\newtheorem{definition}[definition]{Definition}
\aliascntresetthe{definition}

\newaliascnt{remark}{theorem}

\aliascntresetthe{remark}

\usepackage[nameinlink,noabbrev]{cleveref}

\crefname{theorem}{theorem}{theorems}
\Crefname{theorem}{Theorem}{Theorems}

\crefname{lemma}{lemma}{lemmas}
\Crefname{lemma}{Lemma}{Lemmas}

\crefname{proposition}{proposition}{propositions}
\Crefname{proposition}{Proposition}{Propositions}

\crefname{corollary}{corollary}{corollaries}
\Crefname{corollary}{Corollary}{Corollaries}

\crefname{assumption}{assumption}{assumptions}
\Crefname{assumption}{Assumption}{Assumptions}

\crefname{definition}{definition}{definitions}
\Crefname{definition}{Definition}{Definitions}

\crefname{remark}{remark}{remarks}
\Crefname{remark}{Remark}{Remarks}

\crefname{appendix}{appendix}{appendices}
\Crefname{appendix}{Appendix}{Appendices}
\usepackage{pgfplots}
\pgfplotsset{compat=1.18}
\allowdisplaybreaks[1]
\AtBeginDocument{%
  \setlength{\abovedisplayskip}{7pt plus 2pt minus 3pt}%
  \setlength{\belowdisplayskip}{7pt plus 2pt minus 3pt}%
  \setlength{\abovedisplayshortskip}{0pt plus 2pt}%
  \setlength{\belowdisplayshortskip}{5pt plus 2pt minus 2pt}%
}

\hypersetup{pdftitle={Gradient Descent Achieves O(\unichar{949}\unichar{8315}\unichar{8311}\unichar{8260}\unichar{8308}) Convergence},pdfauthor={Chang He and Xiaotian Jiang},pdfsubject={Adaptive stepsizes, restart conditions, and nonconvex gradient descent}}

%% file: marco.tex
\newcommand{\bR}{\mathbb{R}}
\newcommand{\norm}[1]{\lVert #1\rVert}
\newcommand{\inner}[2]{\langle #1,#2\rangle}
\newcommand{\cO}{\mathcal{O}}
\newcommand{\eps}{\varepsilon}

%% file: ref.bib
@book{nesterov2018lectures,
  title={Lectures on Convex Optimization},
  author={Nesterov, Yurii},
  volume={137},
  year={2018},
  publisher={Springer}
}

@inproceedings{carmon2017,
  title={{``Convex Until Proven Guilty'': Dimension-Free Acceleration of Gradient Descent on Non-Convex Functions}},
  author={Carmon, Yair and Duchi, John C. and Hinder, Oliver and Sidford, Aaron},
  booktitle={Proceedings of the 34th International Conference on Machine Learning},
  volume={70},
  series={Proceedings of Machine Learning Research},
  pages={654--663},
  year={2017},
  publisher={PMLR}
}

@article{lilin2023,
  title={Restarted Nonconvex Accelerated Gradient Descent: No More Polylogarithmic Factor in the {$\mathcal{O}(\epsilon^{-7/4})$} Complexity},
  author={Li, Huan and Lin, Zhouchen},
  journal={Journal of Machine Learning Research},
  volume={24},
  number={157},
  pages={1--37},
  year={2023}
}

@inproceedings{agarwal2021fractal,
  title={Acceleration via Fractal Learning Rate Schedules},
  author={Agarwal, Naman and Goel, Surbhi and Zhang, Cyril},
  booktitle={Proceedings of the 38th International Conference on Machine Learning},
  volume={139},
  series={Proceedings of Machine Learning Research},
  pages={87--99},
  year={2021},
  publisher={PMLR}
}

@article{marumo2024,
  author = {Marumo, Naoki and Takeda, Akiko},
  title = {Parameter-Free Accelerated Gradient Descent for Nonconvex Minimization},
  journal = {SIAM Journal on Optimization},
  volume = {34},
  number = {2},
  pages = {2093--2120},
  year = {2024},
  doi = {10.1137/22M1540934}
}
